\pdfoutput=1
\RequirePackage{ifpdf}
\ifpdf % We are running pdfTeX in pdf mode
\documentclass[pdftex]{sigma}
\else
\documentclass{sigma}
\fi

\begin{document}

\allowdisplaybreaks

\renewcommand{\PaperNumber}{043}

\FirstPageHeading

\ShortArticleName{Matrix Weight Function}

\ArticleName{Vector-Valued Polynomials
\\
and a~Matrix Weight Function with $\boldsymbol{B_{2}}$-Action.~II}

\Author{Charles F.~DUNKL}

\AuthorNameForHeading{C.F.~Dunkl}

\Address{Department of Mathematics, University of Virginia,
\\
PO Box 400137, Charlottesville VA 22904-4137, USA}
\Email{\href{mailto:cfd5z@virginia.edu}{cfd5z@virginia.edu}}
\URLaddress{\url{http://people.virginia.edu/~cfd5z/home.html}}

\ArticleDates{Received February 15, 2013, in f\/inal form June 07, 2013; Published online June 12, 2013}

\Abstract{This is a~sequel to [\textit{SIGMA} \textbf{9} (2013), 007, 23 pages], in which there is
a~construction of a~$2\times2$ positive-def\/inite matrix function $K (x )$ on $\mathbb{R}^{2}$.
The entries of $K(x)$ are expressed in terms of hypergeometric functions.
This matrix is used in the formula for a~Gaussian inner product related to the standard module of the
rational Cherednik algebra for the group $W (B_{2} )$ (symmetry group of the square) associated to
the ($2$-dimensional) ref\/lection representation.
The algebra has two parameters: $k_{0}$, $k_{1}$.
In the previous paper $K$~is determined up to a~scalar, namely, the normalization constant.
The conjecture stated there is proven in this note.
An asymptotic formula for a~sum of $_{3}F_{2}$-type is derived and used for the proof.}


\Classification{33C52; 33C20}

\section{Introduction}

This is a~sequel to~\cite{Dunkl2013} and the def\/initions and notations from that paper are used here.
Brief\/ly, we constructed a~$2\times2$ positive-def\/inite matrix function $K(x)$ on
$\mathbb{R}^{2}$ whose entries are expressed in terms of hypergeometric functions.
This matrix is used in the formula for a~Gaussian inner product related to the standard module of the
rational Cherednik algebra for the group $W(B_{2})$ (symmetry group of the square) associated to
the ($2$-dimensional) ref\/lection representation.
The algebra has two parameters: $k_{0}$, $k_{1}$.
In~\cite{Dunkl2013} $K$ is determined up to the normalization constant, henceforth denoted by
$c (k_{0},k_{1} )$.
The conjecture stated there is proven in this note.

Instead of trying to integrate $K$ directly (a problem involving squares of hypergeometric functions whose
argument is $x_{2}^{2}/x_{1}^{2}$) we compute a~sequence of integrals in two ways: asymptotically and
exactly in terms of sums.
Comparing the two answers will determine the value of $c(k_{0},k_{1})$.
First the problem is reduced to a~one-variable integral over the sector
$\left\{ (\cos\theta,\sin\theta ):0<\theta<\frac{\pi}{4}\right\}$ of the unit circle.
With detailed information about the behavior of a~function $f(\theta)$ near $\theta=0$ one can
f\/ind an asymptotic value of $\int_{0}^{\pi/4}\theta^{n}f(\theta)d\theta$.
This part of the argument is described in Section~\ref{integral}.
The other part is produced by exploiting the relationship between the Laplacian and integration over the
circle.
That is, the plan is to determine the result of applying appropriate powers of the Laplacian to certain
polynomials behaving like $\theta^{n}$ near $\theta=0$.
This will be done by establishing recurrence relations; their proofs are in Section~\ref{Thproof}.
The main theorem and its proof which combines the various ingredients are contained in Section~\ref{Eval_c}.

Recall from~\cite[p.~18]{Dunkl2013} that for $0<x_{2}<x_{1}$ and $u=\frac{x_{2}}{x_{1}}$
\begin{gather}
L(u)_{11}=\left\vert u\right\vert^{k_{1}}\left(1-u^{2}\right)^{-k_{0}}F\left(-k_{0},\frac{1}{2}
-k_{0}+k_{1};k_{1}+\frac{1}{2};u^{2}\right),
\label{solL}
\\
L(u)_{12}=-\frac{k_{0}}{k_{1}+\frac{1}{2}}\left\vert u\right\vert^{k_{1}}\left(1-u^{2}
\right)^{-k_{0}}uF\left(1-k_{0},\frac{1}{2}-k_{0}+k_{1};k_{1}+\frac{3}{2};u^{2}\right),
\nonumber
\\
L(u)_{21}=-\frac{k_{0}}{\frac{1}{2}-k_{1}}\left\vert u\right\vert^{-k_{1}}\left(1-u^{2}
\right)^{-k_{0}}uF\left(1-k_{0},\frac{1}{2}-k_{0}-k_{1};\frac{3}{2}-k_{1};u^{2}\right),
\nonumber
\\
L(u)_{22}=\left\vert u\right\vert^{-k_{1}}\left(1-u^{2}\right)^{-k_{0}}F\left(-k_{0},\frac{1}{2}
-k_{0}-k_{1};\frac{1}{2}-k_{1};u^{2}\right),
\nonumber
\end{gather}
and  \cite[p.~20]{Dunkl2013}
\begin{gather}
K(x)=L(u)^{T}\left[
\begin{matrix}d_{1}&0
\\
0&d_{2}
\end{matrix}
\right]L(u),
\label{Kd1d2}
\\
d_{1}=c(k_{0},k_{1})\frac{\Gamma\left(\frac{1}{2}-k_{1}\right)^{2}}{\cos\pi k_{0}
\Gamma\left(\frac{1}{2}+k_{0}-k_{1}\right)\Gamma\left(\frac{1}{2}-k_{0}-k_{1}\right)},
\nonumber
\\
d_{2}=c(k_{0},k_{1})\frac{\Gamma\left(\frac{1}{2}+k_{1}\right)^{2}}{\cos\pi k_{0}
\Gamma\left(\frac{1}{2}+k_{0}+k_{1}\right)\Gamma\left(\frac{1}{2}-k_{0}+k_{1}\right)}.
\nonumber
\end{gather}

\section{Integrals over the circle
\label{integral}
}

First we reduce the integral to a~sector of the unit circle by assuming homogeneity and an invariance
property.
In general, the integral of a~positively homogeneous function with respect to Gaussian measure can be found
by integrating over the sphere (circle).
Furthermore the integral over the sphere of a~polynomial~$p$ homogeneous of degree~$2n$ can be found by
computing~$\Delta^{n}p$.
This method applies in the present situation, when replacing~$\Delta$ by~$\Delta_{\kappa}$.

As before, elements of $\mathcal{P}_{V}$ can be expressed as polynomials
$f_{1}(x)t_{1}+f_{2}(x)t_{2}$ or vectors
$(f_{1}(x),f_{2}(x))$, as needed.
The Gaussian inner product is expressed as
\begin{gather*}
\langle f,g\rangle_{G}=\int_{\mathbb{R}^{2}}f(x)K(x)g(x)^{T}
e^{-\vert x\vert^{2}/2}dx,
\end{gather*}
and the normalization condition is equivalent to
$ \langle(1,0),(1,0) \rangle_{G}=1$.
Recall from~\cite[p.~2]{Dunkl2013} that the fundamental bilinear form $ \langle\cdot,\cdot \rangle_{\tau}$ on
$\mathcal{P}_{V}$ can be written as
\begin{gather*}
 \langle f_{1}(x)t_{1}+f_{2}(x)t_{2},g_{1}(x)t_{1}+g_{2}
(x)t_{2} \rangle_{\tau}
\\
\qquad {} = \langle t_{1},f_{1} (\mathcal{D} ) \{g_{1}(x)t_{1}+g_{2}(x)t_{2}
 \} \rangle_{\tau}|_{x=0}+ \langle t_{2},f_{2} (\mathcal{D} ) \{g_{1}
(x)t_{1}+g_{2}(x)t_{2} \} \rangle_{\tau}|_{x=0},
\end{gather*}
where $f (\mathcal{D} )$ and $|_{x=0}$ denote $f (\mathcal{D}_{1},\mathcal{D}_{2} )$ and
evaluation at $x=0$, respectively.
The abstract Gaussian inner product is then def\/ined as $\langle f,g\rangle_{G}=\langle
e^{\Delta_{\kappa}/2}f,e^{\Delta_{\kappa}/2}g\rangle_{\tau}$ for $f,g\in\mathcal{P}_{V}$.
The key property of this inner product is $ \langle x_{i}f,g \rangle_{G}= \langle
f,x_{i}g \rangle_{G}$ for $i=1,2$.
Subsequently we can derive another abstract inner product which in ef\/fect acts like the use of spherical
polar coordinates in computing integrals with respect to Gaussian measure.
This is done in the following ($S$~symbolizes the sphere/circle):

\begin{definition}
For polynomials $f\in\mathcal{P}_{V,n},g\in\mathcal{P}_{V,m}$ let $\ell=\frac{m+n}{2}$ and
\begin{gather*}
 \langle f,g \rangle_{S}:=\frac{1}{2^{\ell}\ell!} \langle f,g \rangle_{G}=\frac{1}
{2^{\ell}\ell!} \big\langle e^{\Delta_{\kappa}/2}f,e^{\Delta_{\kappa}/2}g \big\rangle_{\tau}
, \qquad m\equiv n \ \operatorname{mod}2,
\\
 \langle f,g \rangle_{S}:=0,\qquad m-n\equiv1\ \operatorname{mod}2.
\end{gather*}
This is extended to all polynomials by linearity.
\end{definition}

Next we relate the Gaussian integral of certain invariant polynomials to the abstract formula for
$ \langle\cdot,\cdot \rangle_{S}$.
Let $x_{\theta}:= (\cos\theta,\sin\theta )$, a~generic point on the unit circle.

\begin{lemma}
\label{polarI}
Suppose $f,g\in\mathcal{P}_{V}$ are relative invariants of the same type, that is, for some linear
character $\chi$ of $W$,
$ (wf )(x)=f (xw )w^{-1}=\chi (w )f(x)$ and
$ (wg )(x)=\chi (w )g(x)$ for each $w\in W$, then
\begin{gather*}
\int_{\mathbb{R}^{2}}f(x)K(x)g(x)^{T}e^{- \vert x \vert^{2}/2}
dx=8\int_{0}^{\infty}e^{-r^{2}/2}rdr\int_{0}^{\pi/4}f (rx_{\theta} )K (x_{\theta}
 )g (rx_{\theta} )^{T}d\theta.
\end{gather*}
\end{lemma}

\begin{proof}
Let $C_{0}= \{x:0<x_{2}<x_{1} \}$, the fundamental chamber, then
\begin{gather*}
\int_{\mathbb{R}^{2}}f(x)K(x)g(x)^{T}e^{- \vert x \vert^{2}/2}
dx=\sum_{w\in W}\int_{C_{0}}f (xw )K (xw )g (xw )^{T}
e^{- \vert x \vert^{2}/2}dx
\\
\hphantom{\int_{\mathbb{R}^{2}}f(x)K(x)g(x)^{T}e^{- \vert x \vert^{2}/2}dx}{}
=\sum_{w\in W}\int_{C_{0}}f (xw  )w^{-1}K(x)wg (xw )^{T}
e^{- \vert x \vert^{2}/2}dx
\\
\hphantom{\int_{\mathbb{R}^{2}}f(x)K(x)g(x)^{T}e^{- \vert x \vert^{2}/2}dx}{}
=\sum_{w\in W}\chi (w )^{2}\int_{C_{0}}f(x)K(x)g(x)^{T}
e^{-  \vert x \vert^{2}/2}dx.
\end{gather*}
The statement follows from the fact $\chi(w)^{2}=1$ and the use of polar coordinates.
Recall $K$ is positively homogeneous of degree zero.
\end{proof}
\begin{proposition}
\label{intproS}
If $f,g\in\mathcal{P}_{V}$ are relative invariants of the same type and
$f\in\mathcal{P}_{V,n},g\in\mathcal{P}_{V,m}$ with $m\equiv n\operatorname{mod}2$ and $\ell=\frac{m+n}{2}$
then
\begin{gather*}
8\int_{0}^{\pi/4}f (x_{\theta} )K (x_{\theta} )g (x_{\theta} )^{T}
d\theta=\frac{1}{2^{\ell}\ell!} \langle f,g \rangle_{G}= \langle f,g \rangle_{S}.
\end{gather*}
If further $n=2q+1$ and $m=1$ then
\begin{gather*}
 \langle f,g \rangle_{S}=\frac{1}{2^{2q+1}q! (q+1 )!} \big\langle\Delta_{\kappa}^{q}
f,g\big\rangle_{\tau}.
\end{gather*}
\end{proposition}

\begin{proof}
From Lemma~\ref{polarI} the factor relating $ \langle\cdot,\cdot \rangle_{S}$ to
$ \langle\cdot,\cdot \rangle_{G}$ is $\int_{0}^{\infty}e^{-r^{2}/2}r^{2\ell+1}dr=2^{\ell}\ell!$.
Now suppose $m=1$ and $\ell=q+1$.
By def\/inition $ \langle f,g \rangle_{S}=\frac{ \langle
e^{\Delta_{\kappa}/2}f,e^{\Delta_{\kappa}/2}g \rangle_{\tau}}{2^{q+1} (q+1 )!}$, then
$e^{\Delta_{\kappa}/2}g=g$ and the degree-$1$ component of $e^{\Delta_{\kappa}/2}f$ is
$\frac{1}{2^{q}q!}\Delta_{\kappa}^{q}f$ (recall that $ \langle f,g \rangle_{\tau}=0$ when $f$ and
$g$ have dif\/ferent degrees of homogeneity).
\end{proof}

We will f\/ind exact expressions in the form of sums for $\Delta_{\kappa}^{q}f$ for certain polynomials in
Section~\ref{Thproof}.

The idea underlying the asymptotic evaluation is this: suppose that $g$ is continuous on $ [0,1 ]$
and satisf\/ies $ \vert g(t)-g(0) \vert\leq Ct$ for some constant, then
$\int_{0}^{1}g(t) (1-t )^{n}dt\allowbreak=\frac{1}{n+1}g(0)+O\left(\frac{1}{n^{2}}\right)$.
This formula can be adapted to the measure $t^{\alpha}dt$ with $\alpha>-1$.
In the sequel we will use~$C$, $C^{\prime}$ to denote constants whose values need not be explicit, as in the
``big $O$'' symbol.

Let $\phi:=x_{1}^{2}-x_{2}^{2}$; then $\phi^{2}$ is $W$-invariant.
Furthermore $\phi^{2n}p_{1,2}$ and $\phi^{2n+1}p_{1,4}$ are all relative invariants of the same type, that
is, $\sigma_{1}f=\sigma_{12}^{+}f=-f$ (recall $p_{1,2}=-x_{2}t_{1}+x_{1}t_{2}$ and
$p_{1,4}=-x_{2}t_{1}-x_{1}t_{2}$).
We will evaluate $ \langle\phi^{2n}p_{1,2},p_{1,2} \rangle_{S}$ and
$ \langle\phi^{2n+1}p_{1,4},p_{1,2} \rangle_{S}$.
These polynomials peak at $\theta=0$ and vanish at $\theta=\frac{\pi}{4}$.
The following are used in the expressions for $p_{1,4}Kp_{1,2}^{T}$ and $p_{1,2}Kp_{1,2}^{T}$.
Set
\begin{gather}
h_{1}(z):=F\left(-k_{0},\frac{1}{2}-k_{0}+k_{1};\frac{3}{2}+k_{1};z\right),
\label{eqhh}
\\
h_{2}(z):=F\left(-k_{0},-\frac{1}{2}-k_{0}-k_{1};\frac{1}{2}-k_{1};z\right),
\nonumber
\\
h_{3}(z):=F\left(k_{0},\frac{1}{2}+k_{0}+k_{1};\frac{3}{2}+k_{1};z\right),
\nonumber
\\
h_{4}(z):=F\left(k_{0},-\frac{1}{2}+k_{0}-k_{1};\frac{1}{2}-k_{1};z\right).
\nonumber
\end{gather}
Each of these hypergeometric series satisf\/ies the criterion for absolute convergence at $z=1$ (for real
$F (a,b;c;z )$ the condition is $c-a-b>0$; here $c-a-b=1\pm2k_{0}$), and so each satisf\/ies
a~bound of the form $ \vert h(z)-1 \vert\leq Cz$ for $0\leq z\leq1$.
Recall the coordinate $u=\frac{x_{2}}{x_{1}}$; on the unit circle (in $-\frac{\pi}{2}<\theta<\frac{\pi}{2}$)
\begin{gather*}
x_{\theta}=\left(\frac{1}{\sqrt{1+u^{2}}},\frac{u}{\sqrt{1+u^{2}}}\right),\qquad \phi=\frac{1-u^{2}}{1+u^{2}}
,\qquad d\theta=\frac{du}{1+u^{2}}.
\end{gather*}
By use of the identities
\begin{gather}
F (a,b;c;z )-\frac{a}{c}zF (a+1,b;c+1,z )=F (a,b-1;c;z ),
\label{FFFz}
\\
F (a,b;c;z )-\frac{a}{c}F (a+1,b;c+1;z )=\frac{c-a}{c}F (a,b;c+1;z ),
\nonumber
\end{gather}
we obtain for $x=x_{\theta}$, $0<\theta<\frac{\pi}{4}$ and $0<u<1$:
\begin{gather}
x_{2}L_{11}-x_{1}L_{12}=\frac{u^{k_{1}+1}\left(1-u^{2}\right)^{-k_{0}}}{\left(1+u^{2}\right)^{1/2}}
\frac{1+2k_{0}+2k_{1}}{1+2k_{1}}h_{1}\left(u^{2}\right),
\label{xLxL}
\\
x_{1}L_{22}-x_{2}L_{21}=\frac{u^{-k_{1}}\left(1-u^{2}\right)^{-k_{0}}}{\left(1+u^{2}\right)^{1/2}}h_{2}
\left(u^{2}\right);
\nonumber
\end{gather}
and
\begin{gather}
-x_{2}L_{11}-x_{1}L_{12}=-\frac{u^{k_{1}+1}\left(1-u^{2}\right)^{k_{0}}}{\left(1+u^{2}\right)^{1/2}}
\frac{1-2k_{0}+2k_{1}}{1+2k_{1}}h_{3}\left(u^{2}\right),
\label{xLp14}
\\
-x_{2}L_{21}-x_{1}L_{22}=-\frac{u^{-k_{1}}\left(1-u^{2}\right)^{k_{0}}}{\left(1+u^{2}\right)^{1/2}}h_{4}
\left(u^{2}\right).
\nonumber
\end{gather}
The expressions for $L_{ij}$ (see~\eqref{solL}) appearing in~\eqref{xLp14} are f\/irst transformed with
$F (a,b;c;z )= (1-z )^{c-a-b}F (c-a,c-b;c;z )$ before using
identities~\eqref{FFFz}.
These formulae will be used to obtain asymptotic expressions for the integrals
$ \langle\phi^{2n}p_{1,2},p_{1,2} \rangle_{S}$ and
$ \langle\phi^{2n+1}p_{1,4},p_{1,2} \rangle_{S}$.
The notation $a (n )\sim b (n )$ means
$\lim\limits_{n\rightarrow\infty}\frac{a (n )}{b (n )}=1$.
\begin{lemma}
Suppose $\alpha,\gamma>-1$, and $n=2,3,\ldots$ then
\begin{gather*}
\int_{0}^{1}t^{\alpha} (1-t )^{n+\gamma} (1+t )^{\beta-n}dt= (2n )^{-\alpha-1}
\Gamma (\alpha+1 )\left(1+O\left(\frac{1}{n}\right)\right).
\end{gather*}
\end{lemma}

\begin{proof}
The term $\left(\frac{1-t}{1+t}\right)^{n}$ is transformed to $\left(1-v\right)^{n}$ by the change of
variable $t=\frac{v}{2-v}$.
The integral becomes
\begin{gather*}
2^{-\alpha-1}\int_{0}^{1}v^{\alpha} (1-v )^{n+\gamma}\left(1-\frac{v}{2}
\right)^{-\alpha-\beta-\gamma-2}dv=2^{-\alpha-1}
\frac{\Gamma (\alpha+1 )\Gamma (n+\gamma+1 )}{\Gamma (n+\alpha+\gamma+2 )}+R,
\end{gather*}
where $R$ is bounded by $2^{-\alpha-1}C\int_{0}^{1}v^{\alpha+1}\left(1-v\right)^{n+\gamma}dv$ and
$\left\vert\left(1-\frac{v}{2}\right)^{-\alpha-\beta-\gamma-2}-1\right\vert\leq Cv$ for $0\leq v\leq1$.
By Stirling's formula $\frac{\Gamma (n+\gamma+1 )}{\Gamma (n+\alpha+\gamma+2 )}\sim
n^{-\alpha-1}$ and $R\sim C^{\prime}n^{-\alpha-2}$ for some constant~$C^{\prime}$.
\end{proof}

\begin{corollary}
Suppose $g(t)$ is continuous and $ \vert g(t)-g(0) \vert\leq Ct$
for $0\leq t\leq1$ then
\begin{gather*}
\int_{0}^{1}t^{\alpha} (1-t )^{n+\gamma} (1+t )^{\beta-n}
g(t)dt= (2n )^{-\alpha-1}
g(0)\Gamma (\alpha+1 )\left(1+O\left(\frac{1}{n}\right)\right).
\end{gather*}
\end{corollary}
\begin{proof}
Break up the integrand as
\begin{gather*}
t^{\alpha} (1-t )^{n+\gamma} (1+t )^{\beta-n}
 \{g(0)+ (g(t)-g(0) ) \};
\end{gather*}
by the Lemma the integral of the second part is bounded by
$C (2n )^{-\alpha-2}\Gamma (\alpha+2 ) \left(1+O\left(\frac{1}{n}\right)\right).$
\end{proof}

\begin{proposition}
\label{intp12}
Suppose $-\frac{1}{2}<k_{0}\pm k_{1}<\frac{1}{2}$ then
\begin{gather*}
\left\langle\phi^{2n}p_{1,2},p_{1,2}\right\rangle_{S}=8\int_{0}^{\pi/4}\phi (x_{\theta} )^{2n}
p_{1,2} (x_{\theta} )K (x_{\theta} )p_{1,2} (x_{\theta} )^{T}d\theta
\\
\hphantom{\left\langle\phi^{2n}p_{1,2},p_{1,2}\right\rangle_{S}}{}
\sim\frac{2\pi c(k_{0},k_{1})}{\cos\pi k_{0}\cos\pi k_{1}}\frac{2^{2k_{1}}n^{k_{1}-1/2}
\Gamma\left(\frac{1}{2}+k_{1}\right)}{\Gamma\left(\frac{1}{2}+k_{0}+k_{1}\right)\Gamma\left(\frac{1}{2}
-k_{0}+k_{1}\right)}.
\end{gather*}
\end{proposition}
\begin{proof}
By def\/inition
\begin{gather*}
\phi (x_{\theta} )^{2n}p_{1,2} (x_{\theta} )K (x_{\theta} )p_{1,2}
 (x_{\theta} )^{T}\\
 \qquad{} =\left(x_{1}^{2}-x_{2}^{2}\right)^{2n}
 \left\{d_{1} (x_{1}L_{12}-x_{2}L_{11} )^{2}+d_{2} (x_{1}L_{22}-x_{2}L_{21} )^{2}
\right\}.
\end{gather*}
Thus equation~\eqref{xLxL} implies
\begin{gather*}
8\int_{0}^{\pi/4}\phi (x_{\theta} )^{2n}p_{1,2} (x_{\theta} )K (x_{\theta}
 )p_{1,2} (x_{\theta} )^{T}d\theta
\\
\qquad{} =8d_{1}\left(\frac{1+2k_{0}+2k_{1}}{1+2k_{1}}\right)^{2}\int_{0}^{1}\left(\frac{1-u^{2}}{1+u^{2}}
\right)^{2n}\frac{\left(1-u^{2}\right)^{-2k_{0}}}{\left(1+u^{2}\right)^{2}}u^{2k_{1}+2}h_{1}\left(u^{2}
\right)^{2}du
\\
\qquad\quad{}
+8d_{2}\int_{0}^{1}\left(\frac{1-u^{2}}{1+u^{2}}\right)^{2n}\frac{\left(1-u^{2}\right)^{-2k_{0}}}
{\left(1+u^{2}\right)^{2}}u^{-2k_{1}}h_{2}\left(u^{2}\right)^{2}du,
\end{gather*}
where $h_{1}$ and $h_{2}$ are from equation~\eqref{eqhh}.
The key fact is that $h_{i}(0)=1$ and $\left\vert h_{i}\left(u^{2}\right)-1\right\vert\leq
Cu^{2}$ for $0\leq u\leq1$ with some constant $C$, $i=1,2$.
In each integral change the variable $u=v^{1/2}$; the f\/irst integral equals
\begin{gather}
\frac{1}{2} (4n )^{-k_{1}-3/2}\Gamma\left(k_{1}+\frac{3}{2}\right)\left(1+O\left(\frac{1}{n}
\right)\right)
\label{intd1}
\end{gather}
and the second integral equals
\begin{gather}
\frac{1}{2} (4n )^{k_{1}-1/2}\Gamma\left(\frac{1}{2}-k_{1}\right)\left(1+O\left(\frac{1}{n}
\right)\right).
\label{intd2}
\end{gather}
This is the dominant term in the sum because $k_{1}-\frac{1}{2}>-k_{1}-\frac{3}{2}$ (that is, $2k_{1}+1>0$).
Using the value of $d_{2}$ in equation~\eqref{Kd1d2} and the identity
$\Gamma\left(\frac{1}{2}-k_{1}\right)\Gamma\left(\frac{1}{2}+k_{1}\right)=\frac{\pi}{\cos\pi k_{1}}$ we
f\/ind
\begin{gather*}
\left\langle\phi^{2n}p_{1,2},p_{1,2}\right\rangle_{S}\sim4c(k_{0},k_{1}) (4n )^{k_{1}
-1/2}\frac{\Gamma\left(\frac{1}{2}-k_{1}\right)\Gamma\left(\frac{1}{2}+k_{1}\right)^{2}}{\cos\pi k_{0}
\Gamma\left(\frac{1}{2}+k_{0}+k_{1}\right)\Gamma\left(\frac{1}{2}-k_{0}+k_{1}\right)}
\\
\hphantom{\left\langle\phi^{2n}p_{1,2},p_{1,2}\right\rangle_{S}}{}
=\frac{2\pi c(k_{0},k_{1})}{\cos\pi k_{0}\cos\pi k_{1}}\frac{2^{2k_{1}}n^{k_{1}-1/2}
\Gamma\left(\frac{1}{2}+k_{1}\right)}{\Gamma\left(\frac{1}{2}+k_{0}+k_{1}\right)\Gamma\left(\frac{1}{2}
-k_{0}+k_{1}\right)}.\tag*{\qed}
\end{gather*}
\renewcommand{\qed}{}
\end{proof}

\begin{proposition}
\label{intp14}
Suppose $-\frac{1}{2}<k_{0}\pm k_{1}<\frac{1}{2}$ then
\begin{gather*}
\left\langle\phi^{2n+1}p_{1,4},p_{1,2}\right\rangle_{S}=8\int_{0}^{\pi/4}\phi (x_{\theta} )^{2n+1}
p_{1,4} (x_{\theta} )K (x_{\theta} )p_{1,2} (x_{\theta} )^{T}d\theta
\\
\hphantom{\left\langle\phi^{2n+1}p_{1,4},p_{1,2}\right\rangle_{S}}{}
\sim\frac{-2\pi c(k_{0},k_{1})}{\cos\pi k_{0}\cos\pi k_{1}}\frac{2^{2k_{1}}n^{k_{1}-1/2}
\Gamma\left(\frac{1}{2}+k_{1}\right)}{\Gamma\left(\frac{1}{2}+k_{0}+k_{1}\right)\Gamma\left(\frac{1}{2}
-k_{0}+k_{1}\right)}.
\end{gather*}
\end{proposition}

\begin{proof}
By def\/inition
\begin{gather*}
\phi (x_{\theta} )^{2n+1}p_{1,4} (x_{\theta} )K (x_{\theta} )p_{1,2}
 (x_{\theta} )^{T}=-\left(x_{1}^{2}-x_{2}^{2}\right)^{2n+1}
\\
\qquad{}\times\left\{d_{1} (x_{2}L_{11}+x_{1}L_{12} ) (x_{2}L_{11}-x_{1}L_{12} )+d_{2} (x_{2}
L_{21}+x_{1}L_{22} ) (x_{2}L_{21}-x_{1}L_{22} )\right\}.
\end{gather*}
Thus equations~\eqref{xLxL} and~\eqref{xLp14} imply
\begin{gather*}
8\int_{0}^{\pi/4}\phi(x_{\theta})^{2n+1}p_{1,4}(x_{\theta})K (x_{\theta})p_{1,2}(x_{\theta})^{T}d\theta
=-8d_{1}\frac{(1-2k_{0}+2k_{1})(1+2k_{0}+2k_{1})}{(1+2k_{1})^{2}}
\\
\qquad{}
\times\int_{0}^{1}\left(\frac{1-u^{2}}{1+u^{2}}\right)^{2n+1}\frac{u^{2k_{1}+2}}{\left(1+u^{2}\right)^{2}}
h_{1}\left(u^{2}\right)h_{3}\left(u^{2}\right)du
\\
\qquad{}
-8d_{2}\int_{0}^{1}\left(\frac{1-u^{2}}{1+u^{2}}\right)^{2n+1}\frac{u^{-2k_{1}}}{\left(1+u^{2}\right)^{2}}
h_{2}\left(u^{2}\right)h_{4}\left(u^{2}\right)du.
\end{gather*}
Arguing as in the previous proof, one obtains the same expressions~\eqref{intd1} and~\eqref{intd2} for the
f\/irst and second integrals respectively.
\end{proof}

\section{Evaluation of the normalizing constant}\label{Eval_c}

The proof of the following appears in Section~\ref{Thproof}.
(Recall $\left(a\right)_{n}:=\prod\limits_{i=1}^{n}\left(a+i-1\right)$.)
\begin{theorem}\label{p14p12}
For arbitrary $k_{0}$, $k_{1}$ and $n\geq0$
\begin{gather*}
\left\langle\phi^{2n}p_{1,2},p_{1,2}\right\rangle_{S}=\frac{1}{n!\left(\frac{1}{2}\right)_{n+1}}
\sum_{j=0}^{n}\frac{(-n)_{j}^{2}}{j!}(-k_{1})_{j}\left(\frac{1}{2}+k_{1}+k_{0}
\right)_{n+1-j}\left(\frac{1}{2}+k_{1}-k_{0}\right)_{n-j},
\\
\left\langle\phi^{2n+1}p_{1,4},p_{1,2}\right\rangle_{S}=-\frac{1}{(n+1)!\left(\frac{1}{2}
\right)_{n+1}}
\\
\hphantom{\left\langle\phi^{2n+1}p_{1,4},p_{1,2}\right\rangle_{S}=}{}
\times\sum_{j=0}^{n}\!\frac{(-n)_{j}(-n{-}1)_{j}}{j!}(-k_{1})_{j}\!\left(\frac{1}
{2}+k_{1}\!+k_{0}\right)_{n+1-j}\! \left(\frac{1}{2}+k_{1}\!-k_{0}\right)_{n+1-j}\!.\!
\end{gather*}
\end{theorem}

\begin{corollary}
\label{p123F2}
Suppose $-k_{1}\pm k_{0}\notin\frac{1}{2}+\mathbb{N}_{0}$ then
\begin{gather*}
\left\langle\phi^{2n}p_{1,2},p_{1,2}\right\rangle_{S}=\frac{\left(\frac{1}{2}+k_{1}+k_{0}\right)_{n+1}
\left(\frac{1}{2}+k_{1}-k_{0}\right)_{n}}{\left(\frac{1}{2}\right)_{n+1}n!}
\\
\hphantom{\left\langle\phi^{2n}p_{1,2},p_{1,2}\right\rangle_{S}=}{}
\times {}_{3}F_{2}\left(\genfrac{}{}{0pt}{}{-n,-n,-k_{1}}{-n-\frac{1}{2}-k_{1}-k_{0},-n+\frac{1}{2}-k_{1}
+k_{0}};1\right),
\\
\left\langle\phi^{2n+1}p_{1,4},p_{1,2}\right\rangle_{S}=-\frac{\left(\frac{1}{2}+k_{1}+k_{0}\right)_{n+1}
\left(\frac{1}{2}+k_{1}-k_{0}\right)_{n+1}}{\left(\frac{1}{2}\right)_{n+1}(n+1)!}
\\
\hphantom{\left\langle\phi^{2n+1}p_{1,4},p_{1,2}\right\rangle_{S}=}{}
\times {}_{3}F_{2}\left(\genfrac{}{}{0pt}{}{-n,-n-1,-k_{1}}{-n-\frac{1}{2}-k_{1}-k_{0},-n-\frac{1}{2}-k_{1}
+k_{0}};1\right).
\end{gather*}
\end{corollary}

The next step is to compare the two hypergeometric series to
$_{2}F_{1}\left(\genfrac{}{}{0pt}{}{-n,-k_{1}}{-n-2k_{1}};1\right)$ which equals
$\frac{(1+k_{1})_{n}}{(1+2k_{1})_{n}}$ (Chu--Vandermonde sum).
The following lemma will be used with $a=\frac{1}{2}+k_{1}+k_{0}$, $b=-\frac{1}{2}+k_{1}-k_{0}$, and
$c=-k_{1}$.

\begin{lemma}\label{sqz}
Suppose $0<a<1$, $-1<b<0$ and $c>-1$ then
\begin{gather*}
_{3}F_{2}\left(\genfrac{}{}{0pt}{}{-n,-n-1,c}{-n-a,-n-b-1};1\right)\leq\frac{\left(1+a+b+c\right)_{n}}
{\left(1+a+b\right)_{n}}\leq {}_{3}F_{2}\left(\genfrac{}{}{0pt}{}{-n,-n,c}{-n-a,-n-b};1\right)\!, \ \ \ c\geq0;
\\
_{3}F_{2}\left(\genfrac{}{}{0pt}{}{-n,-n,c}{-n-a,-n-b};1\right)\leq\frac{\left(1+a+b+c\right)_{n}}
{\left(1+a+b\right)_{n}}\leq {}_{3}F_{2}\left(\genfrac{}{}{0pt}{}{-n,-n-1,c}{-n-a,-n-b-1};1\right)\!,\ \ \ c\leq0.
\end{gather*}
\end{lemma}

\begin{proof}
If $c=0$ then each expression equals $1$.
The middle term equals $_{2}F_{1}\left(\genfrac{}{}{0pt}{}{-n,c}{-n-a-b};1\right)$.
For $0\leq i\leq n$ set
\begin{gather*}
s_{i}:=\frac{(-n)_{i}^{2}(c)_{i}}{i!(-n-a)_{i}(-n-b)_{i}},
\qquad
t_{i}:=\frac{(-n)_{i}(c)_{i}}{i!(-n-a-b)_{i}},
\qquad
u_{i}:=\frac{(-n)_{i}(-n-1)_{i}(c)_{i}}{i!(-n-a)_{i}
(-n-b-1)_{i}}.
\end{gather*}
Note $s_{0}=t_{0}=u_{0}=1$.
From the relation $(-n-d)_{i}=(-1)^{i}(n-i+1+d)_{i}$ it follows that
$\mathrm{sign}(s_{i})=\mathrm{sign}(t_{i})=\mathrm{sign}(u_{i})=\mathrm{sign}((c)_{i})$ for each $i$ with $1\leq i\leq n$ (by hypothesis $a+1>2$ and $b+1>0$).
If $c>0$ then $\mathrm{sign}((c)_{i})=1$ for all $i$ and if $-1<c<0$ then
$\mathrm{sign}((c)_{i})=-1$ for $i\geq1$.
We f\/ind
\begin{gather*}
\frac{s_{i}}{t_{i}}=\frac{(-n)_{i}(-n-a-b)_{i}}{(-n-a)_{i}
(-n-b)_{i}}=\frac{(n-i+1)(n-i+1+a+b)}
{(n-i+1+a)(n-i+1+b)}\frac{s_{i-1}}{t_{i-1}},\qquad i\geq1,
\end{gather*}
and
\begin{gather*}
\frac{m(m+a+b)}{(m+a)(m+b)}=1+\frac{-ab}{(m+a)(m+b)}
>1,\qquad m\geq1,
\end{gather*}
because $-ab>0$ and $b>-1$ (setting $m=n-i+1$).
This shows the sequence $\frac{s_{i}}{t_{i}}$ is positive and increasing.
Also
\begin{gather*}
\frac{u_{i}}{t_{i}}=\frac{(-n-1)_{i}(-n-a-b)_{i}}{(-n-a)_{i}
\left(-n-b-1\right)_{i}}=\frac{\left(n-i+2\right)(n-i+1+a+b)}
{(n-i+1+a)(n-i+2+b)}\frac{u_{i-1}}{t_{i-1}},\qquad i\geq1,
\\
\frac{(m+1)(m+a+b)}{(m+a)(m+b+1)}=1-\frac{b(a-1)}
{(m+a)(m+b+1)}<1,\qquad m\geq1,
\end{gather*}
because $a<1$ and $b<0$ (and $a+b>-1$).
The sequence $\frac{u_{i}}{t_{i}}$ is positive and decreasing.
If $c>0$ then $s_{i},t_{i},u_{i}>0$ for $1\leq i\leq n$ and $\frac{s_{i}}{t_{i}}>1>\frac{u_{i}}{t_{i}}$
implies $s_{i}>t_{i}>u_{i}$.
This proves the f\/irst inequalities.
If $-1<c<0$ then $s_{i},t_{i},u_{i}<0$ for $1\leq i\leq n$ and thus
$\frac{s_{i}}{t_{i}}>1>\frac{u_{i}}{t_{i}}$ implies $s_{i}<t_{i}<u_{i}$.
This proves the second inequalities.
\end{proof}

We will use a~version of Stirling's formula to exploit the lemma:
\begin{gather*}
\frac{\left(a\right)_{n}}{\left(b\right)_{n}}=\frac{\Gamma\left(b\right)}{\Gamma\left(a\right)}
\frac{\Gamma\left(a+n\right)}{\Gamma\left(b+n\right)}\sim\frac{\Gamma\left(b\right)}{\Gamma\left(a\right)}
n^{a-b}.
\end{gather*}
\begin{theorem}
Suppose $-\frac{1}{2}<k_{0}\pm k_{1}<\frac{1}{2}$ then the normalizing constant
\begin{gather*}
c(k_{0},k_{1})=\frac{1}{2\pi}\cos\pi k_{0}\cos\pi k_{1}.
\end{gather*}
\end{theorem}

\begin{proof}
Denote the $_{3}F_{2}$-sums in Corollary~\ref{p123F2} by $f_{1}(n)$ and $f_{2}(n)$
respectively, then by Stirling's formula we obtain
\begin{gather*}
\left\langle\phi^{2n}p_{1,2},p_{1,2}\right\rangle_{S}\sim\frac{\Gamma\left(\frac{1}{2}\right)f_{1}
(n)}{\Gamma\left(\frac{1}{2}+k_{1}+k_{0}\right)\Gamma\left(\frac{1}{2}+k_{1}-k_{0}\right)}
n^{2k_{1}-\frac{1}{2}},
\\
\left\langle\phi^{2n+1}p_{1,4},p_{1,2}\right\rangle_{S}\sim-\frac{\Gamma\left(\frac{1}{2}\right)f_{2}
(n)}{\Gamma\left(\frac{1}{2}+k_{1}+k_{0}\right)\Gamma\left(\frac{1}{2}+k_{1}-k_{0}\right)}
n^{2k_{1}-\frac{1}{2}}.
\end{gather*}
By Propositions~\ref{intp12} and~\ref{intp14} these imply for $i=1,2$
\begin{gather*}
f_{i}(n)\sim\frac{2\pi c(k_{0},k_{1})}{\cos\pi k_{0}\cos\pi k_{1}}\times2^{2k_{1}}
n^{-k_{1}}\frac{\Gamma\left(\frac{1}{2}+k_{1}\right)}{\Gamma\left(\frac{1}{2}\right)}
=\frac{2\pi c(k_{0},k_{1})}{\cos\pi k_{0}\cos\pi k_{1}}n^{-k_{1}}\frac{\Gamma(1+2k_{1})}{\Gamma(1+k_{1})},
\end{gather*}
by the duplication formula.
By Lemma~\ref{sqz} $f_{1}(n)\leq\frac{(1+k_{1})_{n}}{(1+2k_{1})_{n}}\leq
f_{2}(n)$ for $0\leq k_{1}<\frac{1}{2}$, and the reverse inequality holds for
$-\frac{1}{2}<k_{1}<0$.
The fact that
$\frac{(1+k_{1})_{n}}{(1+2k_{1})_{n}}\sim\frac{\Gamma(1+2k_{1})}{\Gamma(1+k_{1})}n^{-k_{1}}$ completes the proof.
\end{proof}

The weight function $K$ is integrable if the inequalities $-\frac{1}{2}<k_{0},k_{1}<\frac{1}{2}$ are
satisf\/ied (and the same value of $c(k_{0},k_{1})$ applies).
However $K$ is not positive-def\/inite and integrable unless $-\frac{1}{2}<k_{0}\pm k_{1}<\frac{1}{2}$.
It was shown in~\cite[p.~21]{Dunkl2013} that $\det K=d_{1}d_{2}$.
By using the (now-known) value of $c(k_{0},k_{1})$ and the values of $d_{1}$ and $d_{2}$
(see~\eqref{Kd1d2}) we f\/ind $\det
K=\frac{1}{4\pi^{2}}\cos\pi\left(k_{0}+k_{1}\right)\cos\pi\left(k_{0}-k_{1}\right)$.

\section[Formulae for $\left\langle\phi^{2n}p_{1,2},p_{1,2}\right\rangle_{S}$ and
$\left\langle\phi^{2n+1}p_{1,4},p_{1,2}\right\rangle_{S}$]{Formulae for $\boldsymbol{\left\langle\phi^{2n}p_{1,2},p_{1,2}\right\rangle_{S}}$ and
$\boldsymbol{\left\langle\phi^{2n+1}p_{1,4},p_{1,2}\right\rangle_{S}}$}\label{Thproof}

The inner products are evaluated by computing $\Delta_{\kappa}^{2n}\left(\phi^{2n}p_{1,2}\right)$ and
$\Delta_{\kappa}^{2n+1}\left(\phi^{2n+1}p_{1,4}\right)$ by means of recurrence relations.
The start is a~product formula for $\Delta_{\kappa}$ (using $\partial_{i}$ to denote
$\frac{\partial}{\partial x_{i}}$):

\begin{lemma}
\label{apDprod}
Suppose $f(x)$ is a~$W$-invariant polynomial and $g\left(x,t\right)\in\mathcal{P}_{V}$ then
\begin{gather}
\Delta_{\kappa} (fg )-f\Delta_{\kappa} (g)=g\Delta f+2 \langle\nabla f,\nabla g \rangle
\label{Dprod1}
\\
\hphantom{\Delta_{\kappa} (fg )-f\Delta_{\kappa} (g)=}{}
+2k_{1}\left(g (x,-t_{1},t_{2} )\frac{\partial_{1}f}{x_{1}}+g (x,t_{1},-t_{2}
 )\frac{\partial_{2}f}{x_{2}}\right)
\nonumber
\\
\hphantom{\Delta_{\kappa} (fg )-f\Delta_{\kappa} (g)=}{}
+2k_{0}\left(g (x,t_{2},t_{1} )\frac{\partial_{1}f-\partial_{2}f}{x_{1}-x_{2}}+g (x,-t_{2}
,-t_{1} )\frac{\partial_{1}f+\partial_{2}f}{x_{1}+x_{2}}\right).
\nonumber
\end{gather}
\end{lemma}

\begin{lemma}\label{deltxsq}
Suppose $f\in\mathcal{P}_{V,2n+1}$ then $\Delta_{\kappa}^{n+1} \vert x \vert^{2}f=4(n+1)(n+2)\Delta_{\kappa}^{n}f$.
\end{lemma}

\begin{proof}
If $g\in\mathcal{P}_{V,m}$ then $\Delta_{\kappa} \vert x \vert^{2}g=4(m+1)g+ \vert
x \vert^{2}\Delta_{\kappa}g$ by~\cite[p.~4, equation~(4)]{Dunkl2013}.
Apply $\Delta_{\kappa}$ repeatedly to this expression, and by induction obtain
$\Delta_{\kappa}^{\ell} \vert
x \vert^{2}g=4\ell (m-\ell+2 )\Delta_{\kappa}^{\ell-1}g+ \vert
x \vert^{2}\Delta_{\kappa}^{\ell}g$ for $\ell=1,2,3,\ldots$.
Set $g=f$, $m=2n+1$ and $\ell=n+1$ then $\Delta_{\kappa}^{n+1}f=0.$
\end{proof}

Recall that $\phi:=x_{1}^{2}-x_{2}^{2}$; and $\phi^{2n}p_{1,2}$ and $\phi^{2n+1}p_{1,4}$ are all relative
invariants of the same type as $p_{1,2}$, that is, $\sigma_{1}f=\sigma_{12}^{+}f=-f$.
Thus $\Delta_{\kappa}^{2n}\left(\phi^{2n}p_{1,2}\right)$ and
$\Delta_{\kappa}^{2n+1}\left(\phi^{2n+1}p_{1,4}\right)$ are both scalar multiples of $p_{1,2}$, because
$\Delta_{\kappa}$ commutes with the action of the group and $p_{1,2}$ is the unique degree-$1$ relative
invariant of this type.
\begin{proposition}
For $n=0,1,2,3,\ldots$
\begin{gather}
\Delta_{\kappa}\phi^{2n}p_{1,2}=-8n(1+2k_{1}+2k_{0})\phi^{2n-1}p_{1,4}
\label{delphi_n}
\\
\hphantom{\Delta_{\kappa}\phi^{2n}p_{1,2}=}{}
+8n(2n-1-2k_{0})\vert x\vert^{2}\phi^{2n-2}p_{1,2},
\nonumber
\\
\Delta_{\kappa}\phi^{2n+1}p_{1,4}=-4(2n+1)(1+2k_{1}-2k_{0})\phi^{2n}p_{1,2}
\label{delph_n2}
\\
\hphantom{\Delta_{\kappa}\phi^{2n+1}p_{1,4}=}{}
+8n(2n+1+2k_{0})\vert x\vert^{2}\phi^{2n-1}p_{1,4}.
\nonumber
\end{gather}
\end{proposition}

\begin{proof}
We use Lemma~\ref{apDprod} with $f=\phi^{2n}$ and $g=p_{1,2}$ or $g=\phi p_{1,4}$.
Simple computation shows that
\begin{gather*}
\Delta\phi^{2n}=8n(2n-1)\vert x\vert^{2}\phi^{2n-2}, \qquad \nabla\phi^{2n}=4n\phi^{2n-1}
\left(x_{1},-x_{2}\right),
\\
\Delta_{\kappa}\left(\phi p_{1,4}\right)=-4(1+2k_{1}-2k_{0})p_{1,2}.
\end{gather*}
For $g=p_{1,2}$ we f\/ind $2\left\langle\nabla\phi^{2n},\nabla g\right\rangle=-8n\phi^{2n-1}p_{1,4}$, the
coef\/f\/icient of $k_{1}$ in~\eqref{Dprod1} is $-16n\phi^{2n-1}p_{1,4}$ and the coef\/f\/icient of $k_{0}$
is $32n\phi^{2n-2}x_{1}x_{2}\left(x_{1}t_{1}-x_{2}t_{2}\right)$.
The f\/irst formula now follows from
\begin{gather*}
x_{1}x_{2}\left(x_{1}t_{1}-x_{2}t_{2}\right)=-\frac{1}{2}\left(\phi p_{1,4}+\vert x\vert^{2}
p_{1,2}\right).
\end{gather*}
For $g=\phi p_{1,4}$ we obtain $2\left\langle\nabla\phi^{2n},\nabla
g\right\rangle=8n\phi^{2n-1}\left(2\vert x\vert^{2}p_{1,4}-\phi p_{1,2}\right)$.
The coef\/f\/icient of $k_{1}$ in~\eqref{Dprod1} is $-16n\phi^{2n}p_{1,2}$ and the coef\/f\/icient of
$k_{0}$ is $-32n\phi^{2n-1}x_{1}x_{2}\left(x_{1}t_{1}+x_{2}t_{2}\right)$.
Similarly to the previous case
\begin{gather*}
x_{1}x_{2}\left(x_{1}t_{1}+x_{2}t_{2}\right)=-\frac{1}{2}\left(\phi p_{1,2}+\vert x\vert^{2}
p_{1,4}\right).
\end{gather*}
The proof of the second formula is completed by adding up the parts, including
$\phi^{2n}\Delta_{\kappa}\phi p_{1,4}$.
\end{proof}

We use this to set up a~recurrence relation.
\begin{definition}
For $n=0,1,2,\ldots$ the constants $\alpha_{n}$, $\beta_{n}$ implicitly depending on $k_{0}$, $k_{1}$ are
def\/ined by
\begin{gather*}
\left\langle\phi^{2n}p_{1,2},p_{1,2}\right\rangle_{S}=\alpha_{n}\left\langle p_{1,2},p_{1,2}
\right\rangle_{S}=\alpha_{n}(1+2k_{1}+2k_{0}),
\\
\left\langle\phi^{2n+1}p_{1,4},p_{1,2}\right\rangle_{S}=\beta_{n}\left\langle p_{1,2},p_{1,2}
\right\rangle_{S}=\beta_{n}(1+2k_{1}+2k_{0}).
\end{gather*}
Also $\alpha_{n}^{\prime}:=2^{4n}(2n)!(2n+1)!\alpha_{n}$ and
$\beta_{n}^{\prime}:=2^{4n+2}(2n+1)!(2n+2)!\beta_{n}$.
\end{definition}

\begin{proposition}
Suppose $n=0,1,2,\ldots$ then
$\Delta_{\kappa}^{2n}\left(\phi^{2n}p_{1,2}\right)=\alpha_{n}^{\prime}p_{1,2}$,
$\Delta_{\kappa}^{2n+1}\left(\phi^{2n+1}p_{1,4}\right)=\beta_{n}^{\prime}p_{1,2}$ and
$\alpha_{n}^{\prime}$, $\beta_{n}^{\prime}$ satisfy the recurrence $($with
$\alpha_{0}^{\prime}=1,\beta_{-1}^{\prime}=0)$
\begin{gather*}
\beta_{n}^{\prime}=-4(2n+1)(1+2k_{1}-2k_{0})\alpha_{n}^{\prime}+64n^{2}
(2n+1)(2n+1+2k_{0})\beta_{n-1}^{\prime},
\\
\alpha_{n}^{\prime}=-8n(1+2k_{1}+2k_{0})\beta_{n-1}^{\prime}+64n^{2}
(2n-1)(2n-1-2k_{0})\alpha_{n-1}^{\prime}, \qquad n\geq1.
\end{gather*}
\end{proposition}

\begin{proof}
By Proposition~\ref{intproS}
\begin{gather*}
2^{4n+1}(2n)!(2n+1)!\left\langle\phi^{2n}p_{1,2},p_{1,2}\right\rangle_{S}\\
\qquad =\left\langle\Delta_{\kappa}^{2n}\left(\phi^{2n}p_{1,2}\right),p_{1,2}\right\rangle_{\tau}
=\alpha_{n} \langle p_{1,2},p_{1,2} \rangle_{\tau}=2\alpha_{n} \langle p_{1,2},p_{1,2}
 \rangle_{S}.
\end{gather*}
Similarly
\begin{gather*}
2^{4n+3}(2n+1)!(2n+2)!\left\langle\phi^{2n+1}p_{1,4},p_{1,2}\right\rangle_{S}\\
\qquad =\left\langle\Delta_{\kappa}^{2n+1}\left(\phi^{2n+1}p_{1,4}\right),p_{1,2}\right\rangle_{\tau}
=\beta_{n} \langle p_{1,2},p_{1,2} \rangle_{\tau}=2\beta_{n} \langle p_{1,2},p_{1,2}
 \rangle_{S}.
\end{gather*}
Apply $\Delta_{\kappa}^{2n-1}$ to both sides of equation~\eqref{delphi_n} to obtain
\begin{gather*}
\alpha_{n}^{\prime}p_{1,2}=-8n(1+2k_{1}+2k_{0})\beta_{n-1}^{\prime}p_{1,2}
+8n(2n-1-2k_{0})\Delta_{\kappa}^{2n-1}\left(\vert x\vert^{2}\phi^{2n-2}p_{1,2}\right).
\end{gather*}
By Lemma~\ref{deltxsq}
\begin{gather*}
\Delta_{\kappa}^{2n-1}\left(\vert x\vert^{2}\phi^{2n-2}p_{1,2}
\right)=8n(2n-1)\Delta_{\kappa}^{2n-2}\left(\phi^{2n-2}p_{1,2}
\right)=8n(2n-1)\alpha_{n-1}^{\prime}p_{1,2}.
\end{gather*}
Apply $\Delta_{\kappa}^{2n}$ to both sides of equation~\eqref{delph_n2} to obtain
\begin{gather*}
\beta_{n}^{\prime}p_{1,2}=-4(2n+1)(1+2k_{1}-2k_{0})\alpha_{n}^{\prime}p_{1,2}
+8n(2n+1+2k_{0})\Delta_{\kappa}^{2n}\left(\vert x\vert^{2}\phi^{2n-1}p_{1,4}\right),
\end{gather*}
and by the same lemma
\begin{gather*}
\Delta_{\kappa}^{2n}\left(\vert x\vert^{2}\phi^{2n-1}p_{1,4}
\right)=8n(2n+1)\Delta_{\kappa}^{2n-1}\left(\phi^{2n-1}p_{1,4}
\right)=8n(2n+1)\beta_{n-1}^{\prime}p_{1,2}.\tag*{\qed}
\end{gather*}
\renewcommand{\qed}{}
\end{proof}

By a~simple computation with the change of scale for $\alpha_{n}$, $\beta_{n}$ we obtain the following:
\begin{corollary}
$\alpha_{n}$, $\beta_{n}$ satisfy the recurrence
\begin{gather*}
\beta_{n}=-\frac{1+2k_{1}-2k_{0}}{2(n+1)}\alpha_{n}+\frac{n(2n+1+2k_{0})}
{(n+1)(2n+1)}\beta_{n-1},
\\
\alpha_{n}=-\frac{1+2k_{1}+2k_{0}}{2n+1}\beta_{n-1}+\frac{2n-1-2k_{0}}{2n+1}\alpha_{n-1}.
\end{gather*}
\end{corollary}

The following formulae arose from examining values of $\alpha_{n}$, $\beta_{n}$ for some small $n$, calculated by
using the recurrence and symbolic computation.
In the following we use relations like $\left( a\right) _{m}\left( a+m\right) =\left( a\right)
_{m+1}$ and $\frac{1}{\left( m-1\right) !}=\frac{m}{m!}$.
\begin{theorem}

Suppose $n=0,1,2,\ldots$ then
\begin{gather*}
\alpha_{n}=\sum_{j=0}^{n}\frac{(-n)_{j}^{2}}{n!\left(\frac{3}{2}\right)_{n}j!}\left(-k_{1}
\right)_{j}\left(\frac{3}{2}+k_{1}+k_{0}\right)_{n-j}\left(\frac{1}{2}+k_{1}-k_{0}\right)_{n-j},
\\
\beta_{n}=-\sum_{j=0}^{n}\frac{(-n)_{j}\left(-1-n\right)_{j}}{(n+1)!\left(\frac{3}{2}
\right)_{n}j!}(-k_{1})_{j}\left(\frac{3}{2}+k_{1}+k_{0}\right)_{n-j}\left(\frac{1}{2}+k_{1}-k_{0}
\right)_{n+1-j}.
\end{gather*}

\end{theorem}

For brevity $k_{+}:=k_{1}+k_{0}$ and $k_{-}:=k_{1}-k_{0}$, as previously.
We use induction on the sequence $\alpha_{0}\rightarrow\beta_{0}\rightarrow
\alpha_{1}\rightarrow\beta_{1}\rightarrow\alpha_{2}\rightarrow\cdots$.
The formulae are clearly valid for $n=0$.
Suppose they are valid for some $n-1$.
Consider $-\frac{1+2k_{1}+2k_{0}}{2n+1}\beta_{n-1}+\frac{2n-1-2k_{0}} {2n+1}\alpha_{n-1}$; split up
the $j$-term in~$\alpha_{n-1}$ by
writing $1=\frac{j-k_{1}}{n-\frac{1}{2}-k_{0}}+\frac{n-\frac{1}{2}+k_{1}-k_{0}
-j}{n-\frac{1}{2}-k_{0}}$, then
\begin{gather}
\frac{2n-1-2k_{0}}{2n+1}\alpha_{n-1}=\frac{1}{\left(n-1\right)!\left(\frac{3}{2}\right)_{n}}\sum_{j=0}
^{n-1}\frac{1}{j!}(1-n)_{j}^{2}\left(\frac{3}{2}+k_{+}\right)_{n-1-j}
\label{anm2a}
\\
\hphantom{\frac{2n-1-2k_{0}}{2n+1}\alpha_{n-1}=}{}
\times\left\{(-k_{1})_{j}\left(\frac{1}{2}+k_{-}\right)_{n-j}+(-k_{1})_{j+1}
\left(\frac{1}{2}+k_{-}\right)_{n-1-j}\right\},
\nonumber
\\
\label{bnm2a}
-\frac{1+2k_{1}+2k_{0}}{2n+1}\beta_{n-1}\\
\qquad{}=\frac{\frac{1}{2}+k_{+}}{n!\left(\frac{3}{2}\right)_{n}}\sum_{j=0}
^{n-1}\frac{1}{j!}(1-n)_{j}(-n)_{j}(-k_{1})_{j}\left(\frac{3}{2}+k_{+}
\right)_{n-1-j}\left(\frac{1}{2}+k_{-}\right)_{n-j}.\nonumber
\end{gather}
The coef\/f\/icient of $ ( -k_{1} ) _{j}$ in the sum of the f\/irst part
of~\{~\} in~\eqref{anm2a} and~\eqref{bnm2a} is
\begin{gather*}
\frac{1}{n!\left(\frac{3}{2}\right)_{n}j!}(1-n)_{j}\left(\frac{3}{2}+k_{+}\right)_{n-1-j}
\left(\frac{1}{2}+k_{-}\right)_{n-j}\left[n(1-n)_{j}+\left(\frac{1}{2}+k_{+}
\right)(-n)_{j}\right]
\\
\qquad{} =\frac{1}{n!\left(\frac{3}{2}\right)_{n}j!}(1-n)_{j}\left(\frac{3}{2}+k_{+}\right)_{n-1-j}
\left(\frac{1}{2}+k_{-}\right)_{n-j}(-n)_{j}\left(n-j+\frac{1}{2}+k_{+}\right)
\\
\qquad{} =\frac{1}{n!\left(\frac{3}{2}\right)_{n}j!}(1-n)_{j}\left(\frac{3}{2}+k_{+}\right)_{n-j}
\left(\frac{1}{2}+k_{-}\right)_{n-j}(-n)_{j}.
\end{gather*}
For $j=0$ this establishes the validity of the $j=0$ term in~$\alpha_{n}$.
For $1\leq j\leq n$ replace~$j$ by $j-1$ in the second part of~\{~\} in~\eqref{anm2a} and obtain
\begin{gather*}
\frac{jn}{n!\left(\frac{3}{2}\right)_{n}j!}(1-n)_{j-1}^{2}\left(\frac{3}{2}+k_{+}\right)_{n-j}
\left(\frac{1}{2}+k_{-}\right)_{n-j}(-k_{1})_{j};
\end{gather*}
adding all up leads to
\begin{gather*}
\frac{1}{n!\left(\frac{3}{2}\right)_{n}j!}\left(\frac{3}{2}+k_{+}\right)_{n-j}\left(\frac{1}{2}+k_{-}
\right)_{n-j}\left\{(1-n)_{j}(-n)_{j}+nj(1-n)_{j-1}^{2}\right\},
\end{gather*}
and the expression in~\{~\} evaluates to $ ( -n) _{j}^{2}$.
This proves the validity of the formula for~$\alpha_{n}$.

To prove the formula for $\beta_{n}$ consider $-\frac{1+2k_{1}-2k_{0} }{2(n+1)
}\alpha_{n}+\frac{n\left( 2n+1+2k_{0}\right) }{(n+1) (2n+1) }\beta_{n-1}$ and as
before split up the $j$-term in $\beta_{n-1}$ by writing $1=\frac{j-k_{1}}{n+\frac{1}
{2}+k_{0}}+\frac{n+\frac{1}{2}+k_{1}+k_{0}-j}{n+\frac{1}{2}+k_{0}}$.
Then
\begin{gather}
-\frac{1+2k_{1}-2k_{0}}{2(n+1)}\alpha_{n}=-\frac{\frac{1}{2}+k_{1}-k_{0}}
{(n+1)!\left(\frac{3}{2}\right)_{n}}\sum_{j=0}^{n}\frac{1}{j!}(-n)_{j}^{2}(-k_{1}
)_{j}\left(\frac{3}{2}+k_{+}\right)_{n-j}\!\left(\frac{1}{2}+k_{-}\right)_{n-j}\!,\!\!
\label{an2b}
\\
\frac{n(2n+1+2k_{0})}{(n+1)(2n+1)}\beta_{n-1}=\frac{1}
{(n+1)!\left(\frac{3}{2}\right)_{n}}\sum_{j=0}^{n-1}\frac{1}{j!}(-n)_{j+1}
(-n)_{j}\left(\frac{1}{2}+k_{-}\right)_{n-j}
\label{bnm2b}
\\
\hphantom{\frac{n(2n+1+2k_{0})}{(n+1)(2n+1)}\beta_{n-1}=}{}
\times\left\{(-k_{1})_{j}\left(\frac{3}{2}+k_{+}\right)_{n-j}+(-k_{1})_{j+1}
\left(\frac{3}{2}+k_{+}\right)_{n-j-1}\right\}.
\nonumber
\end{gather}
The coef\/f\/icient of $ ( -k_{1} ) _{j}$ in the sum of~\eqref{an2b} and the f\/irst part
of~\{~\} in~\eqref{bnm2b} is
\begin{gather*}
\frac{-1}{(n+1)!\left(\frac{3}{2}\right)_{n}j!}(-n)_{j}^{2}\left(\frac{3}{2}+k_{+}
\right)_{n-j}\left(\frac{1}{2}+k_{-}\right)_{n-j}\left(\frac{1}{2}+k_{-}+n-j\right)
\\
\qquad{}=\frac{-1}{(n+1)!\left(\frac{3}{2}\right)_{n}j!}(-n)_{j}^{2}\left(\frac{3}{2}+k_{+}
\right)_{n-j}\left(\frac{1}{2}+k_{-}\right)_{n+1-j}.
\end{gather*}
For $j=0$ this establishes the validity of the $j=0$ term in $\beta_{n}$.
For $1\leq j\leq n$ replace $j$ by $j-1$ in the second part of~\{~\} in~\eqref{bnm2b} and obtain
\begin{gather*}
\frac{j}{(n+1)!\left(\frac{3}{2}\right)_{n}j!}(-n)_{j}(-n)_{j-1}
\left(\frac{1}{2}+k_{-}\right)_{n-j+1}\left(\frac{3}{2}+k_{+}\right)_{n-j};
\end{gather*}
adding all up leads to
\begin{gather*}
\frac{-1}{(n+1)!\left(\frac{3}{2}\right)_{n}j!}(-n)_{j}\left(\frac{1}{2}+k_{-}
\right)_{n-j+1}\left(\frac{3}{2}+k_{+}\right)_{n-j}\left\{(-n)_{j}-j(-n)_{j-1}\right\}
\end{gather*}
and the expression in~\{~\} evaluates to $ ( -1-n ) _{j}$.
This proves the validity of the formula for~$\beta_{n}$ and completes the induction.

This completes the proof of Theorem~\ref{p14p12}.

\pdfbookmark[1]{References}{ref}
\LastPageEnding


\begin{thebibliography}{99}
\footnotesize\itemsep=0pt

\bibitem{Dunkl2013}
Dunkl C.F., Vector-valued polynomials and a~matrix weight function
with $B_{2}$-action, \href{http://dx.doi.org/10.3842/SIGMA.2013.007}{\textit{SIGMA}} \textbf{9} (2013), 007, 23~pages, \href{http://arxiv.org/abs/1210.1177}{arXiv:1210.1177}.

\end{thebibliography}
\end{document}